\documentclass[a4paper,reqno, 10pt]{amsart}
\usepackage[utf8]{inputenc}
\usepackage{tikz}
\usepackage{tikz-cd}
\usepackage[top=2.5cm,bottom=2.5cm,left=2.3cm,right=2.3cm]{geometry}
\usepackage{graphicx}
\usepackage{mathtools}
\usepackage{hyperref}

\usepackage[dvips]{epsfig}
\usepackage{color}
\usepackage{verbatim}
\usepackage{amsgen, amstext,amsbsy,amsopn,amssymb, amsthm, mathptmx, amsfonts,amssymb,amscd,amsmath,nicefrac,euscript,enumerate,url,verbatim,calc,framed}
\usepackage[all]{xy}
\def\red#1{{\color{red}#1}}

\theoremstyle{plain}
\newtheorem{thm}{\bf Theorem}[section]

\newtheorem{prop}[thm]{\bf Proposition}
\newtheorem{lem}[thm]{\bf Lemma}
\newtheorem{cor}[thm]{\bf Corollary}
\newtheorem{setup}{\bf Setup}

\theoremstyle{definition}

\theoremstyle{remark}
\newtheorem{rem}[thm]{\bf Remark}

\newtheorem{note}[thm]{\bf Note}

\DeclareMathOperator{\pd}{proj\,dim}
\DeclareMathOperator{\depth}{depth}
\DeclareMathOperator{\reg}{reg}

\DeclareMathOperator{\pf}{Pf}

\DeclareMathOperator{\Sym}{Sym}

\def\NN{\mathbb{N}}
\def\ZZ{\mathbb{Z}}

\newsavebox\foobox

\begin{document}

\title{\textbf{Almost Complete Intersections: Regularity of powers and Rees algebra}}

\author[Neeraj Kumar]{Neeraj Kumar}
\address{Department of Mathematics, Indian Institute of Technology Hyderabad, Kandi, Sangareddy - 502285, INDIA}
\email{neeraj@math.iith.ac.in, ma19resch11002@iith.ac.in}

\author[Chitra Venugopal]{Chitra Venugopal}

\thanks{{\textit{Date}}: Submitted 16/03/2024, Revised: 21/10/2024}


\begin{abstract}
An almost complete intersection ideal can be seen as a $d$-sequence ideal with the minimal number of generators being one more than its height. In this paper, we give exact formulas for the regularity of powers of graded almost complete intersection ideals satisfying certain conditions. We also present the minimal bigraded free resolutions of Rees algebras associated with linear type ideals having one generator more than their grade and study the properties of Cohen-Macaulayness, Koszulness associated to their diagonals. 
\end{abstract}

\maketitle

\section*{Introduction}

All rings in this paper are considered to be standard graded $K$-algebras ($K$ being a field), unless otherwise specified. Rees algebras and symmetric algebras are important classes of bigraded algebras which have powers and symmetric powers of ideals respectively as its graded components. Corresponding to a graded ideal, the defining equations of the symmetric algebra can be given by the presentation of the ideal, but it is not easy to obtain the associated free resolution.
Whereas in case of Rees algebras, even finding the explicit defining equations are found to be challenging in general. Ideals of linear type are a class of ideals for which the corresponding symmetric algebras and Rees algebras are isomorphic, which is equivalent to the defining equations of Rees algebras being linear in the second set of variables. Complete intersection ideals and $d$-sequence ideals are some classes of ideals of this type. In case of ideals generated by regular sequences, Koszul complexes resolve the ideals and the Eagon-Northcott complex resolves the minimal free resolution of the associated Rees algebra (\cite{EN,Eisenbud}). In \cite{NKCV23}, Rees algebras of certain Pfaffians ideals generated by $d$-sequence which come corresponding to alternating matrices are discussed. There are other special cases for which the corresponding Rees algebras are studied \cite{EN,UV1,RHV,MU,FL,KPU}.

\vspace{2mm}

Another class of ideals of linear type are almost complete intersections. An ideal $I$ in a standard graded $K$-algebra $A$ is said to be an almost complete intersection if the minimal number of generators of $I$ is one more than the height of $I$ along with the property that for all minimal primes $\mathfrak{p}$ of $I$ in $A$, $I_{\mathfrak{p}}$ is a complete intersection. It is proved in \cite{JKS2} that homogeneous almost complete intersection ideals are contained in the class of $d$-sequence ideals (Proposition \ref{aci_d&reg}). Further some study has been done on their regularity \cite{JKS,Shen&Zhu}, combinatorial characterization of graphs which correspond to these ideals (\cite{JKS2}) and their corresponding Rees algebras \cite{HRZ,JKS2}. These classes of ideals can be seen as losing the complete intersection property by one generator, yet unlike the ideals generated by regular sequences, little is known in generality regarding the homological invariants associated to them. 

\vspace{2mm}

Minimal free resolutions of ideals are one of the main tools which helps in studying their structure. Finding these resolutions corresponding to arbitrary ideals in general, is a difficult problem. Applying stronger conditions on the ideal helps in a better understanding of the homological invariants associated to it. In Section \ref{Section3}  (Theorem \ref{aci_resolutions}) we give the minimal bigraded free resolution of the Rees algebra corresponding to an ideal satisfying the following properties: $(1)$ ideal is of linear type, $(2)$ minimal number of generators of the ideal is one more than its grade. This is a generalization of a result in \cite{RBC} extending it to the case of non-equigenerated ideals.   
As a consequence, one obtains resolutions of the bigraded algebras associated to almost complete intersection ideals. In particular, it 
shows that in order to obtain the minimal free resolution of the Rees algebra of almost complete intersection ideals, it suffices to have information about the resolution of $I$.

\vspace{2mm}

A significant algebraic invariant which measures the complexity of a module is the Castelnuovo-Mumford regularity. It is proved in \cite{CHT1999, V2000} that for every homogeneous ideal $I$, the regularity of $I^s$ is of the form $as+b$ for $s \gg 0$ and non-negative constants $a,b$. While computing the value of $b$ is often considered to be a challenging problem, the value of $a$ is widely understood in literature. In Section \ref{s2} (Theorem \ref{powers_aci}), we derive the constant $b$ for powers of almost complete intersections thereby extending the result of Shen and Zhu \cite[Theorem 3.1]{Shen&Zhu} to include the case when the ideal is non-equigenerated.

\vspace{2mm}

Let $c,e \geq 0$ be two integers with $(c,e) \neq (0,0)$. Then the diagonal subalgebra of a bigraded algebra $R$ along the $(c,e)$-diagonal $\Delta= \{(cs,es):s \in \ZZ \}$ is defined as the graded algebra $R_{\Delta}:=\bigoplus_{i \geq 0} R_{(ci,ei)}$ (cf. \cite{STV}). From a geometric point of view, they are the homogeneous
coordinate rings of embeddings of blow-ups of projective varieties along a subvariety. For $c,e \geq 0$, and a homogeneous ideal $I$ in a graded ring $A$, $K[(I^e)_c]$, the $K$-algebra generated by the $c^{th}$ Veronese subring of $I^e$, can be identified as a  $(c,e)$-diagonal of the Rees algebra if $c \geq ed+1$, thereby identifying them as the subalgebras of the Rees algebras in a natural way \cite{STV,CHTV}. This association facilitates an algebraic approach to studying $K[(I^e)_c]$. 

\vspace{2mm}

A standard graded $K$-algebra is Koszul if the minimal free resolution of the corresponding residue field $K$ has a linear resolution. Some articles that have discussed these properties of the diagonal subalgebras of Rees algebras of ideals are \cite{STV,CHTV,CC,Neeraj,AKM} which includes the cases of complete intersections, residual intersections etc. In Section \ref{Section4}, we study the ring-theoretic properties like Cohen-Macaulayness, Koszulness of the diagonals of Rees algebras associated to certain specific classes of ideals of linear type. For a complete intersection ideal $I$ in a polynomial ring minimally generated by $r$ elements of degree $d$, it is proved in \cite[Corollary 6.10]{CHTV} that for $c \geq \frac{d(r-1)}{r}$ and $e>0$, $\mathcal{R}(I)_{\Delta}$ is Koszul and an improved bound of $c \geq \frac{d}{2}$ is given for the case of $r=2$, $d=2$ in \cite[Corollary 3.3]{CC} and for the general case of $r=2$ and any $d$ in \cite[Theorem 3.1]{Neeraj}). In Proposition \ref{koszul}$(1)$, we prove a similar result for an almost complete intersection ideal $I$ minimally generated by $3$ elements in degree $d$ in a graded polynomial ring, which is, $\mathcal{R}(I)_{\Delta}$ is Koszul for $c \geq \frac{d}{3}$ and $e>0$.

\section{Preliminaries} \label{prelim}

The following are some notations used throughout the article.

\begin{enumerate} [a)]
    \item $K$ denotes a field of characteristic zero. 
    \item For an ideal $I$, $\mu(I)$ denotes the minimal number of generators of $I$ and ht$(I)$ denotes the height of $I$.
    \item $A$ denotes a standard graded ring of the form $B/I$ where $B$ is the standard graded polynomial ring over the field $K$ and $I$ is an ideal of $B$.
\end{enumerate}

Let $M$ be a finitely generated graded $A$-module which admits the following graded minimal free resolution

$$ \mathbb{F}: \cdots \stackrel{d_{\mathnormal{i+1}}} \longrightarrow   \bigoplus_{{\mathnormal{j}}} A(-\mathnormal{j})^{{{\beta_{\mathnormal{i},\mathnormal{j}}}}} \stackrel{d_{\mathnormal{i}}} \longrightarrow  \cdots \stackrel{d_{\mathnormal{2}}} \longrightarrow  \bigoplus_{{\mathnormal{j}}} A(-\mathnormal{j})^{{{\beta_{\mathnormal{1},\mathnormal{j}}}}} \stackrel{d_{\mathnormal{1}}}
\longrightarrow \bigoplus_{{\mathnormal{j}}} A(-\mathnormal{j})^{{{\beta_{\mathnormal{0},\mathnormal{j}}}}} \longrightarrow  0.$$

Set  $ t_{\mathnormal{i}}^A(M):=\sup\{ \mathnormal{j} : \beta_{\mathnormal{i},\mathnormal{j}}^A(M)\neq 0\}$ where $\beta_{i,j}$ is the $(i,j)^{\text{th}}$ Betti number of $M$ and by convention, $t_{\mathnormal{i}}^A(M)=-\infty$ if $F_{\mathnormal{i}}=0.$ The Castelnuovo-Mumford regularity is then defined as follows:
$$ \reg_A(M)=sup\{ j-i: \beta_{\mathnormal{i},\mathnormal{j}}^A(M) \neq 0; \, i,j \in \NN\}=sup \{ t_i^A(M)-i:i \in \NN \}.$$

The Castelnuovo-Mumford regularity measures the complexity of the minimal free resolution of a finitely generated $A$-module $M$. We use the notation $\reg(M)$ when the ring $A$ is evident.

Similarly, considering the bigraded minimal free resolution of a finitely generated bigraded module $N$ over a bigraded ring $R$, one obtains biregularities. They are defined as follows,
\begin{enumerate}
\item $\text{reg}_y(N)=\text{sup} \{ b \in \ZZ:\, \beta_{i,(a,b+i)} \neq 0, \text{ for some } i,a \in \ZZ\}$ denotes the $y$-regularity of $N$.
\item $\text{reg}_x(N)=\text{sup} \{ a \in \ZZ:\, \beta_{i,(a+i,b)} \neq 0, \text{ for some } i,b \in \ZZ\}$ denotes the $x$-regularity of $N$.
\end{enumerate}

The following lemma gives an inequality for the regularity of modules in a short exact sequence.

\begin{lem} [Regularity Lemma] \label{reg_lemma} \cite[Lemma 3.1]{Hoa&Tam}
Let $0 \rightarrow M \rightarrow N \rightarrow P \rightarrow 0$ be a short exact sequence of finitely generated graded $A$-modules. Then, the following holds.
\begin{enumerate}
    \item \label{r1} If  $\reg(M) \neq \reg(P)+1,$ then $\reg(N)=\max  \{ \reg(M),\reg(P)\}$.
    \item \label{r2} If  $\reg(N) \neq \reg(P),$ then $\reg(M)=\max  \{ \reg(N),\reg(P)+1\}$.
    \item \label{r3} If  $\reg(M) \neq \reg(N),$ then  $\reg(P)=\max  \{ \reg(M)-1,\reg(N)\}$. 
\end{enumerate}
\end{lem}

Similar to the regularity lemma, one can obtain a lower bound for the depth of a module in terms of that of the free modules in the graded resolution of the same.

\begin{lem}[Depth Lemma] \label{depth_lemma}
Let $M$ be a finitely generated $A$-module, with the free modules in its graded minimal free resolution denoted by $F_i$, $i=0, \ldots, \pd_AM$. Then the following inequality holds,
$$\depth \,(M) \geq \min \{ \depth \,(F_i)-i;\, i \geq 0 \} .$$
\end{lem}

Let $\{\bf{a}\}$$=\{a_1, \ldots, a_n \}$ be a sequence,  $I= \langle a_1, \ldots, a_n  \rangle$ and for $i=1, \ldots, n$,  $I_{i}= \langle a_1, \ldots, a_i \rangle$. Then $\{\bf{a}\}$ is said to be a \textit{$d$-sequence} with respect to a module $M$ if it satisfies the condition $(I_{i-1}M:_Ma_i) \cap IM=I_{i-1}M$ for $i=1, \ldots, n$.

\begin{note} \label{properties}
The following are some properties of $d$-sequences and regular sequences that will be used in this article.
\begin{enumerate}
\item \label{p11} Let $\{f_1, \ldots, f_n\}$ be a $d$-sequence in a ring $A$. Then, $(\langle f_1, \ldots, f_{n-1} \rangle:f_n^s)=(\langle f_1, \ldots, f_{n-1} \rangle:f_n)$ for all $s \geq 1$.
\item \label{p12}  Let $\{f_1, \ldots, f_n\}$ be a $d$-sequence in a ring $A$, and let $I= \langle f_1, \ldots, f_n \rangle$. Then, for $s \geq 1$, $${(\langle f_1, \ldots, f_{i-1} \rangle + I^s:f_i)}=(\langle f_1, \ldots, f_{i-1} \rangle:f_i) + I^{s-1}$$ for all $1 \leq i \leq n$ (\cite[Observation 2.4]{Shen&Zhu}).
\item \label{p2} If $\{f_1, \ldots, f_n\}$ forms a regular sequence, then $(\langle f_1, \ldots, f_{i-1} \rangle:f_i) = \langle f_1, \ldots, f_{i-1} \rangle$ for all $1 \leq i \leq n$.
\end{enumerate}
\end{note}

For a field $K$, $A$ is said to be a standard graded $K$-algebra if $A$ can be identified as $B/I$, where $B$ is a graded polynomial ring over $K$ ($K$-algebra) with the degree of the variable being $1$ and $I$ is its homogeneous ideal. For an equigenerated ideal $\mathfrak{I}$ of $A$, the \textit{Rees algebra of $\mathfrak{I}$ in $A$} is a bigraded $K$-algebra defined as $\mathcal{R}(\mathfrak{I})= \bigoplus_{n\geq0} {\mathfrak{I}}^n$.

\vspace{1mm}

Let $c,e \geq 0$ be two integers with $(c,e) \neq (0,0)$. Then the diagonal subalgebra of a bigraded algebra $R$ along the $(c,e)$-diagonal $\Delta= \{(cs,es):s \in \ZZ \}$ is defined as the graded algebra $R_{\Delta}=\bigoplus_{i \geq 0} R_{(ci,ei)}$ (cf. \cite{STV}).

\begin{rem}  \label{delta-prop}
The following are some results related to the Koszulness and Cohen-Macaulayness of the diagonals of bigraded algebras which will be used in the paper.
\begin{enumerate}
    \item \label{(1)} For $\Delta=(c,e)$ and a bigraded $K$-algebra $R$ with the $i^{th}$ Betti number corresponding to the shift $(a,b)$ being denoted by $\beta_{i,a,b}$, $R_{\Delta}$ is Koszul if $ \max\{ \dfrac{a}{c},\dfrac{b}{e}: \beta_{i,a,b} \neq 0 \} \leq i+1$ for every $i$ (\cite[Theorem 6.2]{CHTV}).
    \item \label{(2)} Let $I= \langle f_1, \ldots,f_s\rangle$ be a homogeneous ideal of a standard graded $K$-polynomial ring $B=K[x_1, \ldots,x_m]$ with deg$(f_j)=d_j$. Let $d=\text{max} \{d_1,\ldots,d_s\}$ and $u= \sum_{i=1}^s d_i$. If $\mathcal{R}(I)$ is Cohen-Macaulay, then for $\alpha=\text{min}\{ (e-1)d+u-m, e(u-m)\}$ and $\beta=\text{min}\{ (e-1)d+u-d_1, e(u-d_1)\}$, $\mathcal{R}(I)_{\Delta}$ is Cohen-Macaulay if $c >\text{max}\{ \alpha,\beta,de\}$ \cite[Theorem 4.5]{Olga}.
\end{enumerate}
\end{rem}

\begin{prop} \cite[Proposition 4.10]{JKS2} \label{aci_d&reg}
Let $I$ be a homogeneous almost complete intersection in a polynomial ring over an infinite field, with the height of $I$ denoted by $n$. Then, $I$ is generated by a homogeneous $d$-sequence $\{f_1, \ldots, f_{n+1}\}$ such that $\{f_1, \ldots, f_{n} \}$ is a regular sequence.    
\end{prop}

\section{Regularity of powers of Almost Complete Intersections} \label{Section2}

In \cite{BHT}, the authors give the regularity of powers of an equigenerated complete intersection $I$ generated by a regular sequence of length $n$ and in degree $d$ as,
\begin{equation*}
\reg(I^s)=ds+(d-1)(n-1).
\end{equation*}

Further, in \cite{Shen&Zhu}, Shen and Zhu give exact formulas for the regularity of powers of equigenerated almost complete intersections.

In this section, we extend their result to give the regularity of powers of almost complete intersections in the non-equigenerated case as well. Our result coincides with \cite[Theorem 3.1]{Shen&Zhu} when the ideal is equigenerated.
Note that an upper bound for the same is discussed in \cite[Corollary 2.11]{JKS}.

\begin{setup} \label{s1}
Let $A$ be a standard graded $K$-algebra and $I$ be a homogeneous ideal of $A$ generated by a $d$-sequence $\{ f_1, \ldots,f_{n}\}$ in degrees $d_1, \ldots, d_n$ with $\{f_1, \ldots,f_{n-1}\}$ being a regular sequence and $J=\langle f_1, \ldots,f_{n-1} \rangle$. Assume that $\reg(A/I) < \sum_{l=1}^{n-1}d_l - n+1$ and $d_n=\max \{ d_i \mid i=1, \ldots,n\}$.
\end{setup}

\begin{lem} \label{lem1}
Consider the notation and assumptions as in Setup \ref{s1}. Then, for $i=0, \ldots,n-1$,
$$ \reg \left( \dfrac{A}{\langle f_1, \ldots,f_i\rangle+I^2} \right)=\sum_{l=1}^{n}d_l - n.$$
\end{lem}

\begin{proof}
Consider the short exact sequence,

$$0 \longrightarrow \dfrac{A}{J:f_{n}}(-d_n) \stackrel {.f_{n}}
    \longrightarrow \dfrac{A}{J} \longrightarrow \dfrac{A}{I} \longrightarrow 0$$

Since $J$ is generated by a regular sequence $\{f_1, \ldots, f_{n-1}\}$ of degrees $d_1, \ldots, d_{n-1}$, respectively, the Koszul complex corresponding to the sequence $\{f_1, \ldots, f_{n-1}\}$ resolves the ideal $J$, and it has the following form:

{\small
\begin{equation} \label{Koszul_res}
0 \longrightarrow A(-(d_{1}+d_{2} + \ldots + d_{n-1})) \stackrel{\phi_{n-1}}{\longrightarrow} \cdots \stackrel{\phi_{k+1}}{\longrightarrow} \bigoplus_{{1 \leq i_1 \leq i_2 \leq \cdots \leq i_k \leq n-1}} A(-(d_{i_1}+d_{i_2} + \ldots + d_{i_k})) \stackrel{\phi_k}{\longrightarrow} \ldots \stackrel{\phi_2}{\longrightarrow} \bigoplus_{j=1}^{n-1} A(-d_j) \stackrel{\phi_1}{\longrightarrow} A \longrightarrow 0,
\end{equation}
}
where $\phi_i's$ are the differentials in the Koszul complex.

Thus, from the resolution  $(\ref{Koszul_res})$, one obtains
$$\reg \left( \dfrac{A}{J} \right) = \sum_{l=1}^{n-1} d_l - (n-1),$$
which is greater than $\reg \left( \dfrac{A}{I} \right)$ by assumption. Hence, 
\begin{equation} \label{eq2}
\reg \left( \dfrac{A}{J:f_{n}} \right) = \sum_{l=1}^{n-1} d_l - (n-1) - d_n.
\end{equation}

The result is now proved by descending induction on $i$. Consider the following short exact sequence for looking at the base case $i=n-1$ of descending induction.

$$0 \longrightarrow \dfrac{A}{(J:f_{n}^2)}(-2d_n) \stackrel {.f_{n}^2}
    \longrightarrow \dfrac{A}{J} \longrightarrow \dfrac{A}{J+\langle f_{n}^2 \rangle} \longrightarrow 0.$$

As a consequence of Note \ref{properties}(\ref{p11}), we have $\reg \left( \dfrac{A}{J:\langle f_{n}^2 \rangle} \right)=\reg \left( \dfrac{A}{J:\langle f_{n} \rangle} \right)=\sum_{l=1}^{n-1}d_l - n+1-d_n.$ Then, 

$\reg \left( \dfrac{A}{J:\langle f_{n}^2 \rangle}(-2d_n) \right)=\reg \left( \dfrac{A}{J} \right) +d_n$.

Therefore, from Lemma \ref{reg_lemma}(\ref{r3}), $\reg \left( \dfrac{A}{J+\langle f_{n}^2 \rangle} \right)=\sum_{l=1}^{n}d_l - n$. Now, assume that $i \leq n-2$ and the statement holds for $i$. Thus, considering the short exact sequence,
$$0 \longrightarrow \dfrac{A}{(\langle f_1, \ldots,f_{i-1}\rangle+I^2):f_i}(-d_i) \stackrel {.f_{i}}
    \longrightarrow \dfrac{A}{\langle f_1, \ldots,f_{i-1}\rangle+I^2} \longrightarrow \dfrac{A}{\langle f_1, \ldots,f_{i}\rangle+I^2} \longrightarrow 0.$$

Since $\{ f_1, \ldots, f_n \}$ forms a $d$-sequence with $\{ f_1, \ldots, f_{n-1} \}$ forming a regular sequence, from Note \ref{properties}(\ref{p12}), we have
$\reg \left( \dfrac{A}{(\langle f_1, \ldots, f_{i-1} \rangle + I^2):f_i} \right)(-d_i) = \reg \left( \dfrac{A}{(\langle f_1, \ldots, f_{i-1} \rangle:f_i) + I} \right)(-d_i).$ From Note \ref{properties}(\ref{p2}), this is equal to $\reg \left( \dfrac{A}{I} \right) + d_i < \sum_{l=1}^{n-1} d_l - (n-1) + d_i \leq \sum_{l=1}^{n} d_l - (n-1) \leq \reg \left( \dfrac{A}{\langle f_1, \ldots, f_{i} \rangle + I^2} \right)+1$. 

Hence, from Lemma \ref{reg_lemma}(\ref{r1}), 

$$\reg \left(\dfrac{A}{\langle f_1, \ldots, f_{i-1} \rangle + I^2} \right) = \reg \left( \dfrac{A}{\langle f_1, \ldots, f_{i} \rangle + I^2} \right) = \sum_{l=1}^{n} d_l - n.$$
\end{proof}

\begin{lem} \label{lem2}
Consider the notation and assumptions as in Setup \ref{s1}. Then, for $s \geq 2$ and $i=0, \ldots,n-1$,
$$ \reg \left( \dfrac{A}{(J:f_n)+\langle f_n^s \rangle} \right)=\sum_{l=1}^{n}d_l-n+(s-2)d_n.$$
\end{lem}

\begin{proof}
First we claim that,
$$ \reg \left( \dfrac{A}{(J:f_n)+\langle f_n^s \rangle} \right)=\reg \left( \dfrac{A}{J:f_n} \right)+sd_n-1.$$

For this consider the following short exact sequence.
 $$0 \longrightarrow \dfrac{A}{J:f_n^{s+1}}(-sd_n) \stackrel {.f_{n}^s}
    \longrightarrow \dfrac{A}{J:f_n} \longrightarrow \dfrac{A}{(J:f_n)+\langle f_n^s \rangle} \longrightarrow 0.$$

From Note \ref{properties}(\ref{p11}), one obtains $\reg \left( \dfrac{A}{J:f_n^{s+1}} \right) (-sd_n)=\reg \left( \dfrac{A}{J:f_n} \right)+sd_n$. Then from Lemma \ref{reg_lemma}(\ref{r3}), 
\begin{equation} \label{eq3}
\reg \left( \dfrac{A}{(J:f_n)+\langle f_n^s \rangle} \right)=\reg \left( \dfrac{A}{J:f_n} \right)+sd_n-1.
\end{equation}

Now substituting the value of $\reg \left( \dfrac{A}{J:f_n} \right)$ from Equation (\ref{eq2}) in Equation (\ref{eq3}), we get the required result.

\end{proof}

From Proposition \ref{aci_d&reg}, it is known that homogeneous almost complete intersection ideals in a polynomial ring over an infinite field have a generating set of $d$-sequence such that removing the last generator results in a regular sequence. Hence, as a consequence of Lemmas \ref{lem1} and \ref{lem2}, one can determine the regularity of powers of almost complete intersections that satisfy certain conditions.

\begin{thm} \label{powers_aci}
   Let $A$ be a standard graded polynomial ring over an infinite field $K$ and $I$ be a homogeneous almost complete intersection ideal of $A$ generated by a $d$-sequence $\{ f_1, \ldots,f_{n}\}$ in degrees $d_1, \ldots, d_n$, where $\{f_1, \ldots,f_{n-1}\}$ is a regular sequence and $J=\langle f_1, \ldots,f_{n-1} \rangle$. Assume that $\reg(A/I) < \sum_{l=1}^{n-1}d_l - n+1$ and that $d_n=\max \{ d_i \mid i=1, \ldots,n\}$
  Then for all $s \geq 2$ and $i=0,\ldots, n-1$, the following is true,
   $$\reg \left( \dfrac{A}{\langle f_1,\ldots,f_{i} \rangle+ I^s}  \right)= \sum_{l=1}^{n}d_l - n +(s-2)d_n.$$
\end{thm}

\begin{proof}
We prove by induction on $s \geq 2$. The base case for $s=2$ is a consequence of Lemma \ref{lem1}. Let $s \geq 3$.
For all $j,s \in \NN $ with $j<n$, we look at the following short exact sequence.
\begin{equation} \label{eq1}
0 \longrightarrow \dfrac{A}{(\langle f_1, \ldots,f_{i}\rangle+I^s):f_{i+1}}(-d_{i+1}) \stackrel {.f_{i+1}}
    \longrightarrow \dfrac{A}{\langle f_1, \ldots,f_{i}\rangle+I^s} \longrightarrow \dfrac{A}{\langle f_1, \ldots,f_{i+1}\rangle+I^s} \longrightarrow 0.
\end{equation}

    We prove the result by descending induction on $i$. Consider the base case of $i=n-1$ of the descending induction in Equation (\ref{eq1}). Then clearly, $\langle f_1, \ldots,f_{(n-1)+1}\rangle+I^s=I$.
    From Note \ref{properties}(\ref{p12}) and Lemma \ref{lem2},
    
    $\reg \left( \dfrac{A}{(J+I^s):f_n}(-d_n) \right) =\reg \left( \dfrac{A}{(J:f_n)+I^{s-1}}\right)+d_n =\reg \left( \dfrac{A}{(J:f_n)+\langle f_n^{s-1} \rangle}\right)+d_n = \sum_{l=1}^{n}d_l - n +(s-3)d_n+d_n=\sum_{l=1}^{n}d_l - n +(s-2)d_n \geq \sum_{l=1}^{n-1}d_l - (n-1)+1 > \reg \left(  \dfrac{A}{I}\right)+1$.
    Hence $\reg \left(\dfrac{A}{J+I^s}\right)=\sum_{l=1}^{n}d_l - n +(s-2)d_n$ by Lemma \ref{reg_lemma}(\ref{r1}).

    Now assume $i \leq n-2$ and that the result holds for $i+1$. In the short exact sequence (\ref{eq1}), by induction hypothesis on $i$, $$\reg \left(\dfrac{A}{\langle f_1, \ldots,f_{i+1}\rangle+I^s}\right)=\sum_{l=1}^{n}d_l - n +(s-2)d_n.$$ Since $\{f_1, \ldots,f_j \}$ forms a $d$-sequence for $1 \leq j \leq n$ and a regular sequence for $1 \leq j \leq n-1$, from Note \ref{properties}(\ref{p12}) and Note \ref{properties}(\ref{p2}), $\reg \left(\dfrac{A}{(\langle f_1, \ldots,f_{i}\rangle+I^s):f_{i+1}}\right)=\reg \left(\dfrac{A}{\langle f_1, \ldots,f_{i}\rangle+I^{s-1}} \right)$. This, by induction hypothesis on $s$, is equal to $\sum_{l=1}^{n}d_l - n +(s-3)d_n$.

    Now, the final result is a consequence of the inequality, $$\reg \left(\dfrac{A}{(\langle f_1, \ldots,f_{i}\rangle+I^s):f_{i+1}}\right)+d_{i+1}\leq \sum_{l=1}^{n}d_l - n +(s-3)d_n+d_n <\reg \left(\dfrac{A}{\langle f_1, \ldots,f_{i+1}\rangle+I^{s}} \right)+1$$ and Lemma \ref{reg_lemma}(\ref{r1}). 
\end{proof}

\begin{note}
The above theorem fails if $d_n$ is not the maximal degree generator of $I$. For example consider $I=\langle f_1,f_2,f_3,f_4\rangle=\langle x_{26}x_{45} - x_{25}x_{46} + x_{24}x_{56}, x_{16}x_{45} - x_{15}x_{46} + x_{14}x_{56}, x_{16}x_{25}x_{34} - x_{15}x_{26}x_{34} - x_{16}x_{24}x_{35} + x_{14}x_{26}x_{35} + x_{15}x_{24}x_{36} - x_{14}x_{25}x_{36},x_{36}x_{45} - x_{35}x_{46} + x_{34}x_{56} \rangle \subset K[x_{14},x_{15},x_{16},x_{24},x_{25},x_{26},x_{34},x_{35},x_{36},x_{45},x_{46},x_{56}]=B$. 

Then $\{f_1,f_2,f_3,f_4\}$ forms a $d$-sequence with $\{f_1,f_2,f_3\}$ being a regular sequence. But $\reg \left( \dfrac{B}{I^3}  \right)= 8$ which is not equal to $
(2+2+2+3)-4+2=7$.    
\end{note}

\begin{note} \label{Note_on_aci}
It is important to note that the formula for the regularity of powers in Theorem \ref{powers_aci} holds true for a general setup of a homogeneous ideal $I$ in a standard graded $K$-algebra $A$, where the field is not necessarily infinite, provided that the following properties are satisfied:
\begin{enumerate}
    \item $I$ has a generating set of $d$-sequence $\{f_1, \ldots, f_n \}$ such that $\{f_1, \ldots, f_{n-1}\}$ is a regular sequence. 
    \item $\reg(A/I) < \sum_{l=1}^{n-1}d_l - n+1$ and $d_n=\max \{ d_i \mid i=1, \ldots,n\}$.
\end{enumerate}
\end{note}

As a consequence of Theorem \ref{powers_aci} and Note \ref{Note_on_aci}, the following can be observed regarding the regularity of powers of almost complete intersections of grade $3$. To understand the set of generators of these ideals, we look at the notion of Pfaffians. 

Given a skew-symmetric matrix $X$, the Pfaffian of $X$ is defined as the square root of the determinant of $X$. For $l \in \NN$, let $\pf_{\bar{l}}(X)$ denote the Pfaffian of $X$ obtained by removing the $l^{th}$ row and the corresponding column of $X$, which can be seen as an element in $B=K[X]$. Let $X$ be skew-symmetric matrix of order $t$ with the following form:
\begin{equation} \label{matrix_X}
X=\begin{bmatrix}
0 & 0 & 0 & x_{14} & x_{15} & \cdots & x_{1t} \\
0 & 0 & 0 & x_{24} & x_{25} & \cdots & x_{2t} \\
0 & 0 & 0 & x_{34} & x_{35} & \cdots & x_{3t} \\
-x_{14} & -x_{24} & -x_{34} &  0 & x_{45} & \cdots & x_{4t} \\
-x_{15} & -x_{25} & -x_{35}  & -x_{45} & 0 & \cdots & x_{5t} \\
\vdots & \vdots & \vdots  & \vdots & \vdots & \cdots & \vdots \\
-x_{1t} & -x_{2t} & -x_{3t}  & -x_{4t} & -x_{5t} & \cdots & 0\\
\end{bmatrix} 
\end{equation}

Then, from the results and discussions in \cite{Brown} and \cite{three_takes_aci}, it can be concluded that an almost complete intersection ideal $I$ in $K[X]$ of grade $3$ and type at least $2$ arises as an ideal generated by Pfaffians of the following form,

\begin{enumerate}
    \item If $t \geq 5$ is an odd integer, then 
    \begin{equation} \label{odd_case}
    I=\langle  \pf_{\bar{1}}(X), \pf_{\bar{2}}(X),\pf_{\bar{3}}(X), \pf_{\overline{123}}(X) \rangle.    
    \end{equation}
    \item If $t \geq 6$ is an even integer, then 
    \begin{equation} \label{even_case}
    I=\langle  \pf(X), \pf_{\overline{12}}(X), \pf_{\overline{13}}(X),\pf_{\overline{23}}(X)\rangle.    
    \end{equation}
\end{enumerate}

\begin{setup} \label{s3}
Let $B=K[X]$ denote a standard graded polynomial ring with $X$ being a skew-symmetric matrix of the form (\ref{matrix_X}). For $t \in \NN$ with $t \geq 5$, let $I$ be an almost complete intersection ideal of $B$ of grade $3$ and type $t-3$ (Equations (\ref{odd_case}) and (\ref{even_case})). 
\end{setup}

\begin{rem} \label{remark_powers}
Note that from \cite[Theorems 1.2 and 1.3]{three_takes_aci}, one can obtain the minimal free resolutions of ideals of the form in Setup \ref{s3}. As a consequence, the following properties can be observed. 
 \begin{enumerate}
\item $B/I$ is Cohen-Macaulay. Since $I$ is perfect in $B$, the observation follows from \cite[Section 19]{Eisenbud}.  
\item \label{reg_main} For $t \geq 5$ and $I$ an ideal of type $t-3$, $\text{reg}(B/I)=t-4$. 
\item For $t \geq 5$ and $I$ an ideal of type $t-3$, $(B/I)^{(c)}$ is Koszul for $c \geq \frac{t}{3}$.
\end{enumerate}
\end{rem}

\begin{lem} \label{d-seq_aci}
Let $X$ be a skew-symmetric matrix of order $t$ of the form given in Equation (\ref{matrix_X}). Then,
\begin{enumerate}
    \item For an odd integer $t \geq 5$, $\{ \pf_{\overline{123}}(X), \pf_{\bar{1}}(X), \pf_{\bar{2}}(X),\pf_{\bar{3}}(X)\}$ is a $d$-sequence such that $\{ \pf_{\overline{123}}(X), \pf_{\bar{1}}(X), \pf_{\bar{2}}(X)\}$ forms a regular sequence in $K[X]$.

    \item For an even integer $t \geq 6$, $\{ \pf_{\overline{12}}(X), \pf_{\overline{13}}(X),\pf_{\overline{23}}(X), \pf(X) \}$ is a $d$-sequence such that $\{ \pf_{\overline{12}}(X), \pf_{\overline{13}}(X),\pf_{\overline{23}}(X)\}$ forms a regular sequence in $K[X]$. 
\end{enumerate}
\end{lem}

\begin{proof}

We have $\pf_{\overline{123}}(X), \pf_{\bar{1}}(X), \pf_{\bar{2}}(X)\}$ forms a regular sequence from \cite[C.2 Lemma]{three_takes_aci}. Further, since $J$ is unmixed, from \cite[Proposition 1.3]{HRZ}, $ \{ \pf_{\overline{123}}(X), \pf_{\bar{1}}(X), \pf_{\bar{2}}(X),\pf_{\bar{3}}(X) \}$ forms a $d$-sequence. 

Similarly $(2)$ can be seen as a consequence of \cite[C.6 Lemma]{three_takes_aci} and \cite[Proposition 1.3]{HRZ}. 
\end{proof}

\begin{cor}
Assume the notation and hypotheses as in Setup \ref{s3}. Then for $m \geq 2$,
\begin{enumerate}
\item If $t \geq 5$ is an odd integer, then for $i=1,2,3$, $$\reg\left (\dfrac{B}{ I^m}\right )=\left(\frac{t-1}{2}\right)(m+2)-5.$$
    
\item If $t \geq 6$ is an even integer, then for $i=1,2,3$, $$\reg \left( \dfrac{B}{I^m} \right)=\left(\frac{t}{2}\right)(m+2)-7.$$
\end{enumerate}
\end{cor}

\begin{proof}
We have $I$ to be an almost complete intersection ideal of $B=K[X]$, a standard graded polynomial ring  with $X$ being a skew-symmetric matrix of the form (\ref{matrix_X}), of grade $3$ and type $t-3$ for $t \in \NN$ and $t \geq 5$. 
\begin{enumerate}
    \item If $t$ is an odd integer, then from Equation (\ref{odd_case}), $I$ is generated by the Pfaffians $ \pf_{\overline{123}}(X), \pf_{\bar{1}}(X), \pf_{\bar{2}}(X),\pf_{\bar{3}}(X)$ which are of degrees $\frac{t-3}{2},\frac{t-1}{2},\frac{t-1}{2},\frac{t-1}{2}$, respectively. From Remark \ref{remark_powers}(\ref{reg_main}), we have $\text{reg}(B/I)=t-4 < \left( \frac{t-3}{2}\right)+ \left( \frac{t-1}{2} \right)+ \left( \frac{t-1}{2} \right)-4+1$. 
    Thus, the result can be seen as a consequence of Lemma \ref{d-seq_aci}, Note \ref{Note_on_aci} and the simplification of the formula derived from Theorem \ref{powers_aci}. 
    \item If $t$ is an even integer, then from Equation (\ref{even_case}), $I$ is generated by the Pfaffians $  \pf_{\overline{12}}(X), \pf_{\overline{13}}(X),\pf_{\overline{23}}(X),  \pf(X)$ which are of degrees $\frac{t-2}{2},\frac{t-2}{2},\frac{t-2}{2},\frac{t}{2}$, respectively. From Remark \ref{remark_powers}(\ref{reg_main}), since $\text{reg}(B/I)=t-4 < \left( \frac{t-2}{2}\right)+ \left( \frac{t-2}{2} \right)+ \left( \frac{t-2}{2} \right)-4+1$, similar to the previous case, the result follows from Lemma \ref{d-seq_aci}, Note \ref{Note_on_aci} and Theorem \ref{powers_aci}. 
\end{enumerate}
\end{proof}

\section{Rees algebra of Almost Complete Intersections} \label{Section3}

This section discusses the resolutions of the Rees algebra of some classes of ideals $I$ of linear type satisfying a relation between their minimal number of generators and their respective heights. 
The resolution is obtained as the mapping cone between the resolution of $I$ and the Buchsbaum-Rim complex, corresponding to the generators of $I$. The idea has been discussed in \cite{RBC} for constructing resolutions of equigenerated ideals. We try to extend the result to include the non-equigenerated case as well.
As a consequence of this, minimal bigarded free resolutions of almost complete intersections are obtained.

The following are the building blocks for constructing the minimal bigraded free resolutions given in Theorem \ref{aci_resolutions}.

Let $I$ be an ideal of a standard graded ring $A$ generated by $\{f_1,\ldots,f_{l}\}$ with $\deg(f_i)=d_i$ and grade $I=l-1$. Let 
$(\mathbb{F.},\phi.)$ denote the graded minimal free resolution of $A/I$. 
Now let $S=A[y_1, \ldots,y_l]$ be a polynomial ring and $F$ and $G$ be free $S$-modules of rank $l$ and $2$ respectively. Then the Buchsbaum-Rim complex (we denote it by $\mathbb{C.}$) corresponding to the matrix $\psi=\begin{bmatrix}
   f_1 & \cdots & f_l \\
   y_1 & \cdots & y_l \\
   \end{bmatrix}$  is defined as follows,

\begin{equation*} \label{BRC}
0 \longrightarrow \wedge^{l}F\otimes (\Sym_{l-3}G)^* \stackrel{\sigma_{l-2}} \longrightarrow \cdots \stackrel{\sigma_2}\longrightarrow \wedge^{3}{F}\otimes (\Sym_{0}G)^* \stackrel{\epsilon}\longrightarrow F \stackrel{\psi} \longrightarrow G.    
\end{equation*}

where 

\begin{enumerate}
    \item $(\Sym_{i}G)^*$ is the dual of the  $i^{th}$ symmetric power of $G$. 
    \item $\epsilon:\wedge^{3}{F}\otimes (\Sym_{0}G)^* \rightarrow F$ is the splice morphism in the Buchsbaum-Rim complex defined as follows. 
    
    Let $ \{ e_1, \ldots, e_{l} \}$ be a basis of $F$ and $\{ g_1,g_2 \}$ be a basis of $G$. We have $(\Sym_{0}G)^*=S$  and for $k= \{ i_1,i_2,i_3; \; i_1 <i_2<i_3\}$, let $e_{i}=e_{i_1} \wedge e_{i_2} \wedge e_{i_3}$ denote an element in $\wedge^{3}{F}$. Then, 
    \begin{equation} \label{epsilon_map}
    \epsilon(e_i)= \sum_{m \subset k, \, |m|=2} \text{sgn}(m \subset k) (\text{det}(\psi_m))e_{k \backslash m} 
    \end{equation}
    where $\text{sgn}(m \subset k)$ is the sign of the permutation of $k$ which puts the elements of $m$ in the first $2$ positions and $\psi_m$ is the minor of order $2$ with the columns indexed by the elements of $m$ \cite{RBC}.
    \item For $2 \leq i \leq l-2$, the maps $\sigma_i$ are defined as follows,
    \begin{equation} \label{sigma_map}
    \sigma_i(\beta \otimes \gamma)=\sum_{s,t}[\psi^*(\beta_s')](\gamma_t') \cdot  \gamma_t'' \otimes \beta_s''.
    \end{equation}
\end{enumerate}

where  $\beta \in (\Sym_{i-1}G)^*$ with $\beta=\sum_s \beta_s' \otimes \beta_s'' \in (G)^* \otimes (\Sym_{i-2}G)^*$ and $\gamma \in \wedge^{i+2}{F}$ with $\gamma=\sum_s \gamma_t' \otimes \gamma_t'' \in \wedge^1 F \otimes \wedge^{i+1}F$.

For $i = 1, \ldots, l$, we have $\psi^*(g_1^*)(e_i) = f_i$ and $\psi^*(g_2^*)(e_i) = y_i$. This clearly indicates that the way the differentials are defined means that the degrees of the generators $f_i$ are crucial in determining the bishifts in the Buchsbaum-Rim complex.

Let, 
\begin{equation} \label{truncated_BRcomplex}
\mathbb{C'.:} \hspace{30mm} \; 0 \longrightarrow \wedge^{l}F\otimes (\text{Sym}_{l-3}G)^* \stackrel{\sigma_{l-2}} \longrightarrow \cdots \stackrel{\sigma_2}\longrightarrow \wedge^{3}{F}\otimes (\text{Sym}_{0}G)^*
\end{equation}
denote the truncated Buchsbaum-Rim complex where $C_i'=\wedge^{i+2}{F}\otimes (\text{Sym}_{i-1}G)^*$, and 
\begin{equation} \label{truncated_res}
\mathbb{F'.:} \hspace{30mm} \ldots \stackrel{\phi_4'} \longrightarrow  \bigoplus_{j \geq 3} S(-j,-1)^{\beta_{3j}} \stackrel{\phi_3'} \longrightarrow  \bigoplus_{j \geq 2} S(-j,-1)^{\beta_{2j}} 
\end{equation}
denote the truncated complex $\mathbb{F.} \otimes S$ where the differentials $\phi_i'$ denotes the base change of $\phi_i$ to $S$.

Next we mention a map of complexes from $\mathbb{C'.}$ to $\mathbb{F'.}$ mentioned in \cite{RBC}.

Let $\nu$ be a map of complexes from the Koszul complex corresponding to $\phi_1$, denoted by $(\wedge^{.} A^l, c.)$, to $(\mathbb{F.},\phi.)$, the minimal bigraded free resolution of $I$, obtained by extending the identity map at degree $0$ and degree $1$ such that the diagram commutes. Further, let $\nu'$ denote the base change of $\nu$ to $S$. This map $\nu'$ helps in defining the required map of complexes $\alpha:  \mathbb{C'.} \rightarrow \mathbb{F'.}$ which has the following form. 

For $1 \leq i \leq l-2$, $\alpha_i: \wedge^{i+2}F \otimes (\Sym_{i-1}G)^* \rightarrow F_{i+1}'$ is given as
the composition of the following two maps.
\begin{equation} \label{moc}
\alpha_i: \wedge^{i+2}F \otimes (\Sym_{i-1}G)^* \stackrel{p \otimes \tau} \longrightarrow \wedge^1 F \otimes \wedge^{i+1}F \stackrel{ \tiny{\begin{bmatrix}
   y_1 & \cdots & y_{l} 
\end{bmatrix}} \otimes \nu_{i+1}'} \longrightarrow S \otimes F_{i+1}' 
\end{equation}
where $p: \wedge^{i+2}F \rightarrow \wedge^1 F \otimes \wedge^{i+1}F$ is the co-multiplication map defined as 
$e_{j_1} \wedge e_{j_2} \wedge \cdots \wedge e_{j_{i+2}} \mapsto \sum_{i} e_{j_i} \otimes e_{j_1} \wedge \ldots \wedge \hat{e_{j_i}} \wedge \ldots \wedge e_{j_{i+2}}$ and $\tau:(\Sym_{i-1}G)^* \rightarrow S$ defined as $g_1^{*(i-1-j)} \otimes g_2^{*(j)} \mapsto \left \{\begin{array}{ccc}
     1        & \text{if} & j=0 \\
     0  & \text{if} & j \geq 1.\\
       \end{array} \right.$$ $

With the notations set, we move on to give the main result in this section.

\begin{thm} \label{aci_resolutions}
Let $B=K[x_1,\ldots, x_m]$ be a standard graded polynomial $K$-algebra and $I=\langle f_1, \ldots, f_{n+1} \rangle$ be a graded ideal of $B$ of linear type with deg $f_j=d_j$ for $j=1, \ldots, n+1$ and grade$(I)=n$ (ht$(I)=n$).
Let $S=K[x_1, \ldots,x_m,y_1, \ldots,y_{n+1}]$ be a bigraded $K$-algebra with deg $x_i=(1,0)$ and deg $y_j=(d_j,1)$ for $i=1, \ldots,m$ and $j=1, \ldots, n+1$.
Then the minimal non-standard bigraded free resolution of $\mathcal{R}(I)$ has the form,
\begin{equation} \label{aci_complex}
 \cdots \stackrel{\delta_4}  \longrightarrow \begin{matrix}         
                       \bigoplus_{j \geq 4} S(-j,-1)^{\beta_{4j}}\\
                       \oplus  \\
                       C_2'
                      \end{matrix}
\stackrel{\delta_3}       \longrightarrow \begin{matrix}
                       \bigoplus_{j \geq 3} S(-j,-1)^{\beta_{3j}}\\
                       \oplus  \\
                       C_1' 
                      \end{matrix} 
\stackrel{\delta_2}       \longrightarrow \begin{matrix}
                       \bigoplus_{j \geq 2} S(-j,-1)^{\beta_{2j}} \\
                      \end{matrix} 
\stackrel{\delta_1}       \longrightarrow S \longrightarrow 0.
\end{equation}

where $\beta_{ij}$ are the $(i,j)^{th}$ Betti numbers of $B/I$ and for $i=1, \dots, n-1$, $$C_i'=\bigoplus_{j=0}^{i-1} \bigoplus_{1 \leq j_1 < j_2< \cdots <j_{i+2} \leq n+1} S(-(d_{j_1}+ d_{j_2}+ \cdots + d_{j_{i+2}}),-j-2)$$  are the modules which occur in the Buchsbaum-Rim complex $\mathbb{C'}.$ corresponding to the matrix $\begin{bmatrix}
   f_1 & \cdots & f_{n+1} \\
   y_1 & \cdots & y_{n+1} \\
   \end{bmatrix}$ with differentials denoted by $\sigma_i$. 
   
Let $\alpha: \mathbb{C'}. \rightarrow  \mathbb{F'}.$ be the map of complexes given in Equation (\ref{moc}). Then for all $i=1, \ldots, \pd_S \mathcal{R}(I)$, the differentials $\delta_i$ are of the following form.

$$ \delta_1= \begin{bmatrix}
             y_1 & \ldots & y_{n+1}
             \end{bmatrix} \circ \phi_2', \; \; \delta_2=\begin{pmatrix}
                      \alpha_1 & \phi_3'
                      \end{pmatrix}, \; \; \delta_i=\begin{pmatrix}
                      \phi_{i+1}' & \alpha_{i-1} \\
                      O & -\sigma_{i-1} 
                      \end{pmatrix}.$$ 

\end{thm}

\begin{proof}
Let $\mathbb{F}.$ be a graded minimal free resolution of $B/I$ of the following form,
\begin{equation} \label{Free_res_of_I}
\cdots \longrightarrow \bigoplus_{j \geq 3}B(-j)^{\beta_{3j}}\stackrel {\phi_3} \longrightarrow \bigoplus_{j \geq 2}B(-j)^{\beta_{2j}} \stackrel{\phi_2}\longrightarrow \bigoplus_{j=1}^{n+1}B(-d_j) \stackrel{\phi_1} \longrightarrow B \longrightarrow 0,  
\end{equation}

and $\mathbb{C}.$ be the Buchsbaum-Rim complex corresponding to the map $\psi:S^{n+1} \rightarrow S^2$ given by the matrix $\begin{bmatrix}
   f_1 & \cdots & f_{n+1} \\
   y_1 & \cdots & y_{n+1} \\
   \end{bmatrix}$ denoted as follows,

\begin{equation} \label{BRC}
0 \longrightarrow \wedge^{n+1}{S^{n+1}}\otimes (\Sym_{n-2}S^2)^* \stackrel{\sigma_{n-1}} \longrightarrow \cdots \stackrel{\sigma_2}\longrightarrow \wedge^{3}{S^{n+1}}\otimes (\Sym_{0}S^2)^* \stackrel{\epsilon}\longrightarrow S^{n+1} \stackrel{\psi} \longrightarrow S^2.    
\end{equation}

Here for $i=1, \ldots,n+1$, $\psi^*(g_1^*)(e_i)=f_i$ and $\psi^*(g_2^*)(e_i)=y_i$ which determines the entries in the matrices corresponding to the differentials $\sigma_i's$.
Then, the way the differentials are defined in Equation (\ref{sigma_map}), one obtains

\begin{equation} \label{BR_shifts}
C_i'=\bigoplus_{j=0}^{i-1} \bigoplus_{1 \leq j_1 < j_2< \cdots <j_{i+2} \leq n+1} S(-(d_{j_1}+ d_{j_2}+ \cdots + d_{j_{i+2}}),-j-2)
\end{equation}
for $i=2, \ldots, n-1$.

Consider $(\wedge^. B^{n+1}, k.)$ to be the Koszul complex corresponding to $\phi_1$.
Now since for $ \{j_1, j_2\} \subseteq [n+1]$ such that $j_1 \neq j_2$, $k_2(e_{j_1} \wedge e_{j_2})=f_{j_1}e_{j_2} -f_{j_2}e_{j_1} \in$ Ker $\phi_1=$ Im $\phi_2$, there exists $\beta_{j_1j_2} \in F_2$ such that $\phi_2(\beta_{j_1j_2})=i(k_2(e_{j_1} \wedge e_{j_2}))$ where $i$ is the natural map from $\wedge^{0}B^{n+1}$ to $F_1$. Let $\alpha:\mathbb{C'.} \rightarrow \mathbb{F'.}$ be the map of complexes given in Equation (\ref{moc}) where $\mathbb{C'.}$ and $\mathbb{F'.}$ are the complexes given in Equations (\ref{truncated_BRcomplex}) and (\ref{truncated_res}) respectively. Then $\alpha_1:\wedge^{3}{S^{n+1}}\otimes (\Sym_{0}S^2)^* \rightarrow F_2'$ will be defined as, for $\{j_1,j_2,j_3 \} \subseteq [n+1]$, $$\alpha_1((e_{j_1} \wedge e_{j_2} \wedge e_{j_3})\otimes 1)=y_{j_1}\beta_{j_2j_3}-y_{j_2}\beta_{j_1j_3}+y_{j_3}\beta_{j_1j_2}.$$

Thus, from the way the map $\alpha_1$ is defined, one obtains,
\begin{equation}\label{BR_shift1}
C_1' = \bigoplus_{1 \leq j_1 < j_2 < j_3 \leq n+1} S(-(d_{j_1}+d_{j_2}+d_{j_3}),-2).
\end{equation}

With the details above, one obtains the required resolution (\ref{aci_complex}) by considering the mapping cone of $\alpha$ and following arguments similar to those in \cite{RBC}. 
From the construction of the Complex (\ref{aci_complex}), it is evident that the bishifts in the resolution of the Rees algebra of $I$ is derived from the information in the Complexes (\ref{Free_res_of_I}) and (\ref{BRC}) and Equations (\ref{BR_shifts}) and (\ref{BR_shift1}). 
Further, the resolution is minimal since the entries in the matrices given by the differentials $\delta_i$ are from $\mathfrak{m}S$, where $\mathfrak{m}$ is the bigraded maximal ideal of $S$. Note that the linear shifts in the second variable are contributed by the differentials in the minimal free resolution of $B/I$, while the higher degree shifts are contributed by the $\alpha_i's$ and $\sigma_i's$.  
\end{proof}

\begin{note}
Let $\beta_i$ be the $i^{th}$ total Betti number of $B / I$. Then from the way the maps are defined it is clear that $d_i$, $\alpha_i$, $\delta_i$ represents matrices of order $\beta_{i-1} \times \beta_{i}$, $\beta_{i+1} \times i \cdot {{n+1}\choose{i+2}}$ and $(i-1) \cdot {{n+1}\choose{i+1}} \times i \cdot {{n+1}\choose{i+2}}$ respectively.
\end{note}

\begin{rem} \label{remark}
Note that in Theorem \ref{aci_resolutions}
if $I$ is considered to be a graded ideal of $B$ with $\mu(I)=\text{ht}(I)+1$ and not necessarily of linear type, then the Complex (\ref{aci_complex}) gives a minimal bigraded free resolution of the $\Sym(I)$.
\end{rem}

\begin{rem} \label{pd}
Let $I$ be an graded ideal of $B$ of linear type such that $\mu(I)=\text{ht}(I)+1$. Then, 
$$\pd_S \mathcal{R}(I)=\left \{\begin{array}{ccc}
     n         & \text{if} & \text{pd } (B/I)=n \\
     k-1  & \text{if} & \text{pd } (B/I)=k, \, k > n.\\
       \end{array} \right.$$
\end{rem}


The following observations can be made on the nature of the generators of $I$ as a consequence of the explicit minimal bigraded free resolutions of the symmetric algebra of $I$ obtained from Theorem \ref{aci_resolutions}.

Note that for an ideal $I$ in $B$, its generating set $\{a_1, \ldots, a_n\}$ forms a \textit{strong $s$-sequence} (w.r.t an admissible term order for the monomials in $y_i$ with $y_1 <y_2 <\ldots <y_n$) if for $\Sym(I) \cong S/J$ where $S$ is a bigraded polynomial ring $S=B[y_1,\ldots,y_n]$ and $J$ the defining ideal of the symmetric algebra, $in(J)= \langle L_1y_1,  \ldots, L_ny_n  \rangle $ where $L_i=(a_1,\ldots,a_{i-1}:a_i)$ and $L_1 \subseteq L_2 \subseteq \ldots \subseteq L_n$ (cf. \cite{s-seq}).

\begin{prop}
Let $I$ be graded ideal of a standard graded polynomial ring $B$ satisfying the assumptions in Remark \ref{remark}. Then,
\begin{enumerate}
\item $I$ is generated by a strong $s$-sequence.
\item $I$ is generated by a $d$-sequence if $I$ is of linear type.
\end{enumerate}
\end{prop}

\begin{proof}
From the minimal bigraded free resolution of the $\Sym(I)$ in Complex (\ref{aci_complex}), it is clear that the $y$-regularity of the $\Sym(I)$ is $0$. The result then follows from \cite[Corollary 3.2(i)]{Rom}.

Further when the ideal is of linear type, from Complex (\ref{aci_complex}), $y$-regularity of the $\mathcal{R}(I)$ vanishes. Then $I$ is generated by a $d$-sequence follows from \cite[Corollary 3.2(ii)]{Rom}.
\end{proof}

Some important classes of ideals for which \red{Proposition \ref {aci_resolutions}} gives the bigraded free resolutions of the corresponding symmetric algebras are the following. 

\begin{enumerate}
    \item If $I$ is a graded ideal of $B$ of the form $I=J+\langle g \rangle$ where $J$ is a complete intersection and $g \in B$ such that $I$ is minimally generated by the generators of the ideal $J$ and $g$.
    \item $d$-sequence ideals satisfying the condition that the grade of the ideal is one more than the minimal number of generators of the ideal. In particular, this class includes almost complete intersections.
   
\end{enumerate}

In the remaining part of this section, we give the explicit resolutions of the Rees algebra of some classes of ideals  satisfying the hypothesis of Theorem \ref{aci_resolutions}, in particular almost complete intersections of grade $2$, grade $3$ etc.

\subsection{Almost complete intersections of grade $2$ } \label{subsection1}
\text{ }

\begin{setup} \label{s2}
Let $B$ be a standard graded polynomial ring and $I=\langle f_1, f_2,f_3 \rangle \subseteq B$ be an almost complete intersection of grade $2$ and perfect with deg $f_i=d_i$, $i=1,2,3$. 
\end{setup}

\begin{prop} \label{aci_grade2}
   Consider the ideal $I$ as in Setup \ref{s2} and the notations as in Theorem \ref{aci_resolutions}. Let $\beta_{ij}$ denote the degree of the $(i,j)^{th}$ element in the presentation matrix of $I$. Then the minimal bigraded free resolution of the corresponding $\mathcal{R}(I)$ is of the following form:

  $$ 0 \longrightarrow    \begin{matrix}
                          S(-(d_1+d_2+d_3),-2)
                          \end{matrix} 
   \stackrel{\delta_2} \longrightarrow    \begin{matrix}
                          S(-(d_k+\beta_{k1}),-1)\\
                          \oplus  \\
                          S(-(d_l+\beta_{l2}),-1)
                          \end{matrix}                     \stackrel{\delta_1} 
                          \longrightarrow    S \longrightarrow 0.
$$
where $\beta_{k1}$ and $\beta_{l2}$ are non-zero.
\end{prop}

\begin{proof}
From the Hilbert-Burch theorem \cite[Theorem 20.15]{Eisenbud}, the graded minimal free resolution of $B/I$ has the form,
 $$0 \longrightarrow \begin{matrix}
                          B(-(d_k+\beta_{k1}))\\
                          \oplus  \\
                          B(-(d_l+\beta_{l2}))
                          \end{matrix} 
    \stackrel {\phi_2}   \longrightarrow    \begin{matrix}
                          B(-d_1) \\
                          \oplus  \\
                          B(-d_2) \\
                          \oplus  \\
                          B(-d_3) 
                          \end{matrix} 
                          \stackrel {\phi_1} \longrightarrow  B \longrightarrow 0.$$

where $\phi_1= \left[\begin{matrix}
               f_1 & f_2 & f_3
               \end{matrix}\right]$
               and $\phi_2= \left[ \begin{matrix}
                      z_{11} & z_{12} \\
                      z_{21} & z_{22} \\
                      z_{31} & z_{32} \\
                      \end{matrix} \right]$ and $\beta_{k1}$ and $\beta_{l2}$ come corresponding to some non-zero entries $z_{k1}$ and $z_{l2}$ respectively in $\phi_2$.
       
With a fixed ordered basis, the Buchsbaum-Rim complex ($\mathbb{C}$-complex) has the form, 
$$0 \longrightarrow \begin{matrix}
                          S(-(d_1+d_2+d_3),-2)\\
                          \end{matrix} 
    \stackrel {\epsilon} \longrightarrow  S^{3} \stackrel {\psi}\longrightarrow S^2 $$ where $\psi= \begin{bmatrix}
   f_1 & f_2 & f_3 \\
   y_1 & y_2 & y_3 \\ 
   \end{bmatrix}$ and $ \epsilon=\begin{bmatrix}
                               f_2y_3-f_3y_2 \\
                               f_3y_1-f_1y_3 \\
                               f_1y_2-f_2y_1 
                               \end{bmatrix}.$

Let $z$ be a non-zero divisor in $B$ such that $f_1=-z(z_{21}z_{32}-z_{22}z_{31}),\, f_2=z(z_{11}z_{32}-z_{12}z_{31}), \text{ and } f_3=-z(z_{11}z_{22}-z_{12}z_{21})$, which exists from \cite{Eisenbud},  then $\alpha_1= \begin{bmatrix}
                               -z(z_{12}y_1+z_{22}y_2+z_{32}y_3) \\
                               z(z_{11}y_1+z_{21}y_2+z_{31}y_3)
                                 \end{bmatrix}.$

The result thus follows from Theorem \ref{aci_resolutions}.
\end{proof}

\begin{rem} \label{grade2_rmk}
Perfect ideals of grade $2$ are characterized by the Hilbert-Burch theorem \cite{Eisenbud}. Thus, any $3 \times 2$ matrix $X$ with entries from a graded polynomial ring $B$ that satisfies the condition that the grade of the ideal of the corresponding maximal minors $I_2(X)$ is bounded below by $2 $, along with a non-zero divisor $z$ in $B$, corresponds to an ideal $zI_2(X)$, which has the form mentioned in Setup \ref{s2}. Consequently, the associated Rees algebra is resolved using Proposition \ref{aci_grade2}. 
\end{rem}

\subsection{Almost complete intersections of grade $3$} \label{subsection2} \text{ }

\vspace{1mm}

\begin{prop}
Consider the ideal $I$ as in Setup \ref{s3}. Then the minimal bigraded free resolution of $\mathcal{R}(I)\cong S/J $ where $S=K[X,Y]$ with $Y=\begin{bmatrix}
   y_1, \ldots, y_4 \\
   \end{bmatrix}$ and deg $x_{kl}=(1,0)$ and deg $y_i=(d_i,1)$, $i=1, \ldots, 4$ and $J$ is the defining ideal of $\mathcal{R}(I)$, has the following form:

\begin{enumerate}
    \item When $t=2r+1$, $r \in \NN$, $r \geq 2$, 
    $$ 0 \longrightarrow \begin{matrix}
                          S(-4r+1,-2)\\
                          \oplus  \\
                          S(-4r+1,-3)
                          \end{matrix} 
    \stackrel {\delta_3}  \longrightarrow    \begin{matrix}
                          S(-2r,-1)^{t-3} \\
                          \oplus  \\
                          S(-3r+1,-2)^{3} \\
                          \oplus  \\
                          S(-3r,-2) 
                          \end{matrix} 
                          \stackrel {\delta_2} \longrightarrow    \begin{matrix}
                          S(-2r+1,-1)^t    \\
                          \end{matrix}                       
       \stackrel{\delta_1} \longrightarrow    S \longrightarrow 0.
$$
\item When $t=2r$, $r \in \NN$, $r \geq 3$, 
    $$ 0 \longrightarrow \begin{matrix}
                          S(-4r+3,-2)\\
                          \oplus  \\
                          S(-4r+3,-3)
                          \end{matrix}
    \stackrel {\delta_3'}
       \longrightarrow    \begin{matrix}
                          S(-2r+1,-1)^{t-3} \\
                          \oplus  \\
                          S(-3r+3,-2)^{3} \\
                          \oplus  \\
                          S(-3r+2,-2)^{} 
                          \end{matrix} 
    \stackrel {\delta_2'}
       \longrightarrow    \begin{matrix}
                          S(-2r+2,-1)^t    \\
                          \end{matrix}    
       \stackrel{\delta_1'}           
       \longrightarrow    S \longrightarrow 0.
$$
\end{enumerate}

\end{prop}

\begin{proof}

Let $I$ be an almost complete intersection of grade $3$ and type $t-3$. The following can be observed from the results in \cite{three_takes_aci}.

\begin{enumerate}
    \item If $t=2r+1$, $r \geq 2$, then the graded minimal free resolution of $B/I$ is given by 
     $$ 0 \longrightarrow \begin{matrix}
                          B(-2r)^{t-3}
                          \end{matrix} 
      \stackrel{\phi_3} \longrightarrow    \begin{matrix}
                          B(-2r+1)^{t} 
                          \end{matrix} 
      \stackrel{\phi_2} \longrightarrow    \begin{matrix}
                          B(-r)^3\\
                          \oplus   \\
                          B(-r+1)
                          \end{matrix}                       
      \stackrel{\phi_1} \longrightarrow    B \longrightarrow 0. $$

\item If $t=2r$, $r \geq 3$, then the graded minimal free resolution of $B/I$ is given by 
     $$ 0 \longrightarrow \begin{matrix}
                          B(-2r+1)^{t-3}
                          \end{matrix} 
      \stackrel{\phi_3} \longrightarrow    \begin{matrix}
                          B(-2r+2)^{t
                          } 
                          \end{matrix} 
      \stackrel{\phi_2} \longrightarrow    \begin{matrix}
                          B(-r)    \\
                          \oplus   \\
                          B(-r+1)^3
                          \end{matrix}                       
      \stackrel{\phi_1} \longrightarrow    B \longrightarrow 0. $$
\end{enumerate}

Let $f_{ij}$ denote the element $f_iy_j-f_jy_i$ in $S$. With a fixed ordered basis, the respective $\mathbb{C.}$-complexes have the form,

\begin{enumerate}
\item $$0 \longrightarrow \begin{matrix}
                         S(-4r+1,-2) \\
                          \oplus     \\
                          S(-4r+1,-3) \\
                          \end{matrix}
    \stackrel {\sigma_2} \longrightarrow \begin{matrix}
                          S(-3r+1,-2)^3\\
                          \oplus     \\
                          S(-3r,-2) \\
                          \end{matrix}
    \stackrel {\epsilon} \longrightarrow  S^{3} \stackrel {\psi}\longrightarrow S^2 $$  
\item $$0 \longrightarrow \begin{matrix}
                         S(-4r+3,-2) \\
                          \oplus     \\
                          S(-4r+3,-3) \\
                          \end{matrix}
    \stackrel {\sigma_2} \longrightarrow \begin{matrix}
                          S(-3r+3,-2)^3\\
                          \oplus     \\
                          S(-3r+2,-2) \\
                          \end{matrix}
    \stackrel {\epsilon} \longrightarrow  S^{3} \stackrel {\psi}\longrightarrow S^2 $$
\end{enumerate}
where $\psi= \begin{bmatrix}
   f_1 & f_2 & f_3 & f_4 \\
   y_1 & y_2 & y_3 & y_4\\ 
   \end{bmatrix}$, $\epsilon=\begin{bmatrix}
                               f_{23} & f_{24}  & f_{34} &    0   \\
                              -f_{13} & -f_{14} &    0   & f_{34} \\
                               f_{12} &  0      &-f_{14} & -f_{24}\\
                                  0   & f_{12}  & f_{13} & f_{23} \\
                               \end{bmatrix}$ and $\sigma_2=\begin{bmatrix}
                                        -f_4 & -y_4 \\
                                         f_3 & y_3 \\
                                        -f_2 & -y_2 \\
                                         f_1 & y_1 \\
                                        \end{bmatrix}.$
\vspace{2mm}

Then the bigraded free resolutions of $\mathcal{R}(I)$ are consequences of Theorem \ref{aci_resolutions}.

\end{proof}

\begin{cor}
The bi-regularities of $\mathcal{R}(I)$ when $I$ is an almost complete intersection of grade $3$ of type $t-3$ in a standard graded polynomial ring are as follows,
    \begin{enumerate}
        \item reg$_y\mathcal{R}(I)=0$.
        \item reg$_x\mathcal{R}(I)=4r-4$ when $t=2r+1$, $r \geq 2$.
        \item reg$_x\mathcal{R}(I)=4r-6$ when $t=2r$, $r \geq 3$. 
    \end{enumerate}
\end{cor}

\subsection{Restriction on the degree of the generators of the ideal and certain colon ideals} \label{subsection3} \text{ }
Let $I \subseteq B$ such that ht$(I)=n$ and $\mu(I)=n+1$. Let $I= \langle f_1, \ldots, f_{n+1} \rangle$ such that $f_1, \ldots, f_{n}$ be $B$-regular and let $J= \langle f_1, \ldots, f_{n} \rangle$ with deg $f_i=d$ for all $i$. \\

The $\mathbb{C.}$-complex corresponding to $\psi=\begin{bmatrix}
      f_1 & \cdots & f_{n+1} \\
      y_1 & \cdots & y_{n+1} \\
     \end{bmatrix}$ has the form,

$$ 0  \longrightarrow \bigoplus_{l=2}^n S(-(n+1)d,-l) \stackrel{\sigma_{n-1}} \longrightarrow \cdots \stackrel{\sigma_3} \longrightarrow \begin{matrix}
                          S(-4d,-2)^{n+1 \choose 4}\\
                          \oplus  \\
                          S(-4d,-3)^{n+1 \choose 4}\\
                          \end{matrix} 
    \stackrel {\sigma_2}   \longrightarrow    \begin{matrix}
                          S(-3d,-2)^{n+1 \choose 3}\\
                          \end{matrix} 
                          \stackrel {\epsilon} \longrightarrow
                          S^{n+1}       \stackrel{{\psi}} \longrightarrow S^2.$$

Now to get a free resolution of $B/I$, consider the following short exact sequence, 
\begin{equation} \label{ses1}
0 \longrightarrow \dfrac{B}{(J:f_{n+1})}(-d) \stackrel {.f_{n+1}}
    \longrightarrow \dfrac{B}{J} \longrightarrow \dfrac{B}{I} \longrightarrow 0.
\end{equation}

Since $J$ is generated by a regular sequence, the corresponding Koszul complex resolves $B/J$, that is, the graded minimal free resolution of $B/J$ is given by,

$$ 0 \longrightarrow B(-nd) \stackrel {\eta_{n}}
    \longrightarrow B(-(n-1)d)^{n} \stackrel {\eta_{n-1}} \longrightarrow \cdots \longrightarrow B(-2d)^{n\choose2} \stackrel {\eta_{2}} \longrightarrow B(-d)^n \stackrel {\eta_{1}} \longrightarrow B \longrightarrow 0. 
$$

\begin{enumerate}[(A)]
    \item \label{a} Assume $J'=(J:f_{n+1})$ to be an \textit{equigenerated perfect ideal of grade $2$} with the degree of the generators being $d'$, $\mu(J')=t$, $t \in \NN$ and $J'$ being linearly presented. By the Hilbert-Burch theorem, the graded minimal free resolution of $B/J'$ has the following form:
    $$ 0 \longrightarrow B(-(d'+1))^{t-1} \stackrel {\zeta_{2}}
    \longrightarrow B(-d')^{t}\stackrel {\zeta_{1}} \longrightarrow B \longrightarrow 0 $$

Then the mapping cone of the map, multiplication by $f_{n+1}$ in the short exact sequence (\ref{ses1}), gives a graded free resolution of $B/I$ (not necessarily minimal) of the following form:

\small{$$ 0 \longrightarrow B(-nd) \stackrel  {\phi_{n}}
    \longrightarrow B(-(n-1)d)^{n} \stackrel {\phi_{n-1}} \longrightarrow \cdots \stackrel {\phi_{4}}\longrightarrow  \begin{matrix}
            B(-(d'+d+1))^{t-1}\\
            \oplus \\
            B(-3d)^{{n\choose{3}}}  \\
            \end{matrix}
            \stackrel {\phi_{3}} \longrightarrow \begin{matrix}
            B(-(d'+d))^t \\
            \oplus \\
            B(-2d)^{{n\choose2}}  \\
    \end{matrix} 
     \stackrel {\phi_{2}} \longrightarrow \begin{matrix}
      B(-d)^{n+1} \\
     \end{matrix} \stackrel {\phi_{1}} \longrightarrow B \longrightarrow 0 $$}

\normalsize

As a consequence of Theorem \ref{aci_resolutions}, the following would be a (non-standard) bigraded free resolution of $\mathcal{R}(I)$,

\footnotesize{
$$ 0 \longrightarrow \bigoplus_{l=2}^n S(-(n+1)d,-l) \stackrel{\delta_{n}}
\longrightarrow  \begin{matrix} 
                      S(-nd,-1)  \\
                      \oplus     \\
                     \bigoplus_{l=2}^{n-1} S(-nd,-l)^{n+1}
                      \end{matrix}
               \stackrel {\delta_{n-1}} \longrightarrow \cdots  
    \stackrel {\delta_3}   \longrightarrow    \begin{matrix}
                          S(-(d'+d+1),-1)^{t-1} \\
                          \oplus \\
                          S(-3d,-1)^{{n \choose 3}} \\
                          \oplus  \\
                          S(-3d,-2)^{n+1 \choose 3}\\
                          
                          \end{matrix} 
                          \stackrel {\delta_2} \longrightarrow    \begin{matrix}
                          S(-(d'+d),-1)^t \\
                          \oplus \\
                          S(-2d,-1)^{{n \choose 2}}   \\
                          \end{matrix}                       
       \stackrel{{\delta_1}} \longrightarrow    S \longrightarrow 0.
$$}

\normalsize

\begin{rem}

\begin{enumerate} [1)]
    \item The ideal $J'$ need not always be equigenerated. For example let $I=\langle f_1, f_2, f_3 \rangle$, an ideal in a polynomial ring $B$, be a complete intersection ideal. Then the defining relations of the corresponding Rees algebra is given by the maximal minors of the matrix $X=\begin{bmatrix}
    f_1&f_2&f_3 \\
    y_1&y_2&y_3
\end{bmatrix}$. Let the relations be denoted by $\{g_1,g_2,g_3\}$. Then they satisfy the condition that $\{g_1,g_2\}$ are $S$-regular and $(\langle g_1,g_2 \rangle:g_3)=\langle f_3,y_3\rangle$ (\cite[Lemma 3.1]{CC}), which is clearly not equigenerated.

\item An example of a class of ideals $J'$ that satisfies the condition (\ref{a}) are complete intersections minimally generated by $2$ linear forms. For example, for $I=\langle x_1x_2,x_3x_4,x_2x_3\rangle \subset K[x_1,x_2,x_3,x_4]$, the set $\{x_1x_2,x_3x_4\}$ forms a regular sequence, and $(\langle x_1x_2,x_3x_4\rangle: x_2x_3)= \langle x_1,x_4  \rangle$ is a grade $2$ perfect ideal. In general, as mentioned in Remark \ref{grade2_rmk}, all perfect ideals of grade $2$ are characterized by the Hilbert Burch theorem.

\end{enumerate}
\end{rem}

    \item \label{b} Let  $J'=(J:f_{n+1})$ be an \textit{equigenerated Gorenstein ideal of grade $3$} with the degree of the generators being $d'$, $\mu(J')=t$ and the ideal $J'$ being linearly presented. Then the graded minimal free resolution of $B/J'$ has the following form:
    $$ 0 \longrightarrow B(-(2d'+1)) \stackrel {\zeta_{1}^*}\longrightarrow B(-(d'+1))^{t} \stackrel {\zeta_{2}}
    \longrightarrow B(-d')^{t}\stackrel {\zeta_{1}} \longrightarrow B \longrightarrow 0 $$

As mentioned in (\ref{a}), by considering the mapping cone of the map, multiplication by $f_{n+1}$ in the short exact sequence (\ref{ses1}), one obtains a graded free resolution of $B/I$ (not necessarily minimal) of the following form:
\small{$$ 0 \longrightarrow B(-nd) \stackrel  {\phi_{n}}
    \longrightarrow  \cdots \longrightarrow \begin{matrix}
            B(-(2d'+d+1)) \\
            \oplus \\
            B(-4d)^{{n\choose4}}  \\
            \end{matrix} \stackrel {\phi_{4}}\longrightarrow  \begin{matrix}
            B(-(d'+d+1))^{t} \\
            \oplus \\
            B(-3d)^{{n\choose3}}  \\
            \end{matrix}
            \stackrel {\phi_{3}} \longrightarrow \begin{matrix}
            B(-(d'+d))^t  \\
            \oplus \\
            B(-2d)^{{n\choose2}}  \\
    \end{matrix} 
     \stackrel {\phi_{2}} \longrightarrow \begin{matrix}
      B(-d)^{n+1} \\
     \end{matrix} \stackrel {\phi_{1}} \longrightarrow B \longrightarrow 0 $$}

\normalsize 
From Theorem \ref{aci_resolutions}, the following would be a bigraded free resolution of $\mathcal{R}(I)$,

\footnotesize{
$$ 0  \longrightarrow \bigoplus_{l=2}^n S(-(n+1)d,-l) \stackrel{\delta_{n}}
\longrightarrow  \cdots \longrightarrow \begin{matrix}
                          S(-(2d'+d+1),-1) \\
                          \oplus \\
                          S(-4d,-1)^{{n \choose 4}} \\
                          \oplus \\
                          S(-4d,-2)^{n+1 \choose 4}\\
                          \oplus  \\
                          S(-4d,-3)^{n+1 \choose 4}\\
                          \end{matrix} 
    \stackrel {\delta_3}   \longrightarrow    \begin{matrix}
                          S(-(d'+d+1),-1)^t \\ 
                          \oplus \\
                          S(-3d,-1)^{{n \choose 3}} \\
                          \oplus  \\
                          S(-3d,-2)^{n+1 \choose 3}\\
                          \end{matrix} 
                          \stackrel {\delta_2} \longrightarrow    \begin{matrix}
                          S(-(d'+d),-1)^t \\
                          \oplus \\
                          S(-2d,-1)^{{n \choose 2}}   \\
                          \end{matrix}                       
       \stackrel{{\delta_1}} \longrightarrow    S \longrightarrow 0.
$$}

\end{enumerate}

\begin{note}
In both cases (\ref{a}) and (\ref{b}), we have $d' \leq d$ since $J \subseteq J'$. This should be kept in mind while studying the properties of diagonals in the next section.
\end{note}

\normalsize

\section{Application to diagonal subalgebras of almost complete intersections} \label{Section4}

By extending the idea of the Segre product of algebras, Simis, Trung, and Valla introduced the concept of diagonal subalgebras in \cite{STV}. Finding appropriate conditions on a bigraded algebra $R$ such that some algebraic features of $R$ are inherited by $R_{\Delta}$ is one of the main difficulties in the study of diagonal subalgebras.

Priddy introduced Koszul algebras in 1970 (\cite{Priddy}), in the study of homological properties of graded algebras derived from various constructions in algebraic topology. One of the reasons that makes the study of Koszul algebras important is their frequent appearances in classical constructions in commutative algebra and algebraic geometry. For a Koszul algebra, the coresponding Poincar\'e series is found to be rational (cf. \cite{Froberg}). Moreover, Avramov, Eisenbud, and Peeva gave the characterization, $R$ is Koszul $\Longleftrightarrow$ reg$_R(K) < \infty$ $\Longleftrightarrow$ reg$_R(M) < \infty$ for every finitely generated $R$-module $M$ (cf. \cite{AE,AP}), analogous to the result of Auslander, Buchsbaum, and Serre characterising regular rings in terms of the finiteness of the projective dimension of modules over them. A good survey on Koszul algebras is given in \cite{Froberg,CNR}. 

\vspace{2mm}

Consider the notations as in Theorem \ref{aci_resolutions}.

\subsection{Koszulness of $\mathcal{R}(I)_{\Delta}$}

Let $I= \langle f_1, \ldots, f_{n+1} \rangle$ be an equigenerated ideal of $B$ satisfying the conditions in Theorem \ref{aci_resolutions} with deg $f_i=d$. Then it is possible to associate a standard bigraded structure to the corresponding Rees algebra by considering it to be the quotient of a standard bigraded polynomial ring $S$ with deg $x_i=(1,0)$ and deg $y_j=(0,1)$ for $i=1, \ldots,m$ and $j=1, \ldots, n+1$. From the form of the differentials mentioned in Theorem \ref{aci_resolutions}, the standard bigraded free resolution of the Rees algebra will take the following form,

\begin{equation} \label{aci_standardcomplex}
 \cdots \stackrel{\delta_4}  \longrightarrow \begin{matrix}         
                       \bigoplus_{j \geq 3} S(-j+d,-1)^{\beta_{4j}}\\
                       \oplus  \\
                       C_2' 
                      \end{matrix}
\stackrel{\delta_3}       \longrightarrow \begin{matrix}
                       \bigoplus_{j \geq 2} S(-j+d,-1)^{\beta_{3j}}\\
                       \oplus  \\
                       C_1' 
                      \end{matrix} 
\stackrel{\delta_2}       \longrightarrow \begin{matrix}
                       \bigoplus_{j \geq 1} S(-j+d,-1)^{\beta_{2j}} \\
                      \end{matrix} 
\stackrel{\delta_1}       \longrightarrow S \longrightarrow 0.
\end{equation}
where $$C_i'=\bigoplus_{j=0}^{i-1} S(-(i-j)d,-j-2)^{{n+1}\choose{i+2}}$$ for $i=1, \ldots, n-1$.
Having defined standard bigraded structure on $\mathcal{R}(I)$, it makes sense to talk about the Koszulness of the corresponding diagonals, since $\Delta$ being an exact functor, when applied to the standard bigraded resolution of $\mathcal{R}(I)$ gives a standard graded resolution of $\mathcal{R}(I)_{\Delta}$. 

Let $i$ denote the homological degree and
$t_{i}^B(B/I)$ be as defined in Section \ref{prelim}. Let $ a_{\max}^i$ and $b_{\max}^i$ denote the maximum shifts in the first and second degrees, respectively, at the $i^{th}$ homological degree of the Complex (\ref{aci_standardcomplex}). Clearly, 
$$
a_{\max}^{i} = \max\{ t_{i+1}^B(B/I) + d, (i-1)d \}
$$
and 
$$
b_{\max}^{i} = \max\{j + 2 \mid j = 0, \ldots, i-1\} = i + 1.
$$

Then from Remark \ref{delta-prop}(\ref{(1)}), $\mathcal{R}(I)_{\Delta}$ is Koszul if $\dfrac{a_{max}^{i}}{c} \leq i+1$ and $\dfrac{b_{\max}^{i}}{e} \leq i + 1$  for all $i \geq 1$.

Since $b_{\max}^{i}=i+1$, the condition $\dfrac{b_{\max}^{i}}{e} \leq i + 1$ 
holds true for all $e > 0$. It remains to check the values of $c$ for which 
$
\dfrac{a_{\max}^{i}}{c} \leq i + 1.
$
For this, we need to consider the following two cases from the Complex (\ref{aci_standardcomplex}),

\noindent \textbf{Case 1:} $a_{max}^{i}=(i-1)d$. Then $\dfrac{a_{max}^{i}}{c} \leq i+1 \Longleftrightarrow \dfrac{(i-1)d}{c} \leq i+1 \Longleftrightarrow \dfrac{i-1}{i+1}d \leq c \Longleftrightarrow c \geq d $ since $\left\lvert{\dfrac{i-1}{i+1}}\right \rvert\leq 1$. Thus for $c \geq d$ and $e>0$, $\mathcal{R}(I)_{\Delta}$ is Koszul.

\vspace{2mm}

\noindent \textbf{Case 2:} $a_{max}^{i}=t_{i+1}^B(B/I)-d$. Then for $c \geq  \dfrac{t_{i+1}^B(B/I)-d}{i+1}$ and $e>0$, $\mathcal{R}(I)_{\Delta}$ is Koszul. 

\noindent Let $\gamma =\text{max}\{\dfrac{t_{i+1}^B(B/I)-d}{i+1} | \, 1 \leq i \leq \pd_B(B/I) \}$. Then for $c \geq \max \{ \gamma, d\} $ and $e>0$, $\mathcal{R}(I)_{\Delta}$ is Koszul. 

\vspace{2mm}

As a consequence of the discussion above, the following can be observed regarding the Koszulness of the diagonals of $\mathcal{R}(I)$ for the class of ideals discussed in Section \ref{Section3}.

\begin{prop} \label{koszul}
Let $I \subseteq B$ be a graded almost complete intersection ideal of $B$ and $S$ be a standard bigraded polynomial $K$-algebra. 
\begin{enumerate}
    \item Consider the assumptions as in Setup \ref{s2} with $d_1=d_2=d_3=d$, $d>1$ and $I$ being linearly presented. Then $\mathcal{R}(I)_{\Delta}$ is Koszul for $c \geq \frac{d}{3}$ and $e > 0$.
    \item Let $I$ be an equigenerated ideal of height $n$ with a generating set $\{f_1, \ldots, f_{n+1}\}$ such that $f_1, \ldots, f_n$ forms a regular sequence. Further, let $J=\langle f_1, \ldots, f_n \rangle$
    and $J'=(J:f_{n+1})$ satisfy the assumptions as in  Subsection \ref{subsection3} (\ref{a}) or (\ref{b}). Then for $c \geq \big( \frac{n-1}{n}\big)d$ and $e >0$, $\mathcal{R}(I)_{\Delta}$ is Koszul. 
    
\end{enumerate}
\end{prop}

\begin{rem}
If $I$ is an almost complete intersection ideal of grade $2$ and perfect generated in degree $1$, then $\mathcal{R}(I)_{\Delta}$ is Koszul for all $\Delta$.  
\end{rem}

Note that Proposition \ref{koszul} provides bounds on the diagonals of the Rees algebra associated with almost complete intersections, for which they are Koszul, similar to those found for the Rees algebra corresponding to complete intersections \cite{CHTV,CC,Neeraj}, residual intersections \cite{AKM} and maximal order Pfaffians \cite{NKCV23}.

\subsection{Cohen-Macaulayness of $\mathcal{R}(I)_{\Delta}$}

The following result gives the diagonals for which the diagonal subalgebras of the Rees algberas corresponding to ideals discussed in Section \ref{Section3} would be Cohen-Macaulay.

\begin{prop}
Let $I=\langle f_1, \ldots, f_{n+1}\rangle$ be a graded ideal in a standard graded polynomial ring $B=K[x_1, \ldots, x_m]$ with the degree of $f_i=d_i$, satisfying either of the following conditions:
\begin{enumerate}
\item $I$ be an almost complete intersection ideal satisfying the assumptions as in Setup \ref{s3} or Setup \ref{s2}.

\item The height of the ideal $I$ is $n$ with $f_1, \ldots, f_{n}$ forming a regular sequence. Further, let $J=\langle f_1, \ldots, f_n \rangle$ and $J'=(J:f_{n+1})$ satisfy the assumptions as in  Subsection \ref{subsection3} (\ref{a}) or (\ref{b}).
\end{enumerate}
Set $d=\max\{d_1, \ldots, d_{n+1}\}$ and $u=\sum_{i=1}^{n+1}d_i$. Then $\mathcal{R}(I)$ is Cohen-Macaulay and $\mathcal{R}(I)_{\Delta}$ is Cohen-Macaulay for all $e>0$ and $c >\max \{ \min\{ (e-1)d+u-m, e(u-m)\}, \min\{ (e-1)d+u-d_1, e(u-d_1)\}, ed\}$.  
\end{prop}

\begin{proof}
We have $\dim S=m+n+1=\depth S$ and $\dim \mathcal{R}(I)=m+1$. Since $\pd_B(B/I)=n$ in all the cases, from Remark \ref{pd} and Lemma \ref{depth_lemma}, one obtains $\depth \mathcal{R}(I) \geq m+1$. Hence $\mathcal{R}(I)$ is Cohen-Macaulay. Then, Cohen-Macaulayness of $\mathcal{R}(I)_{\Delta}$ for the mentioned diagonals can be seen as a consequence of Remark \ref{delta-prop}(\ref{(2)}).  
\end{proof}

\vspace{2mm}

\noindent {\bf Acknowledgement.} The Core Research Grant (CRG/2023/007668) from the Science and Engineering Research Board, ANRF, India, partially supports the first author. INSPIRE Fellowship, DST, India, financially supports the second author.

\end{document}